\input amstex
\documentstyle{amsppt}
\magnification=1200
\pageheight{18.5cm}
\pagewidth{12.5cm}

\NoBlackBoxes
\TagsOnRight

\define\cc#1{\Bbb C^{#1}}
\define\ch#1{\Bbb {CH}^{#1}}
\define\hq#1{\Bbb {HQ}^{#1}}
\define\ir#1{\Bbb R^{#1}}
\define\is#1{\Bbb S^{#1}}
\define\bb#1{\Bbb B^{#1}}
\define\grs#1#2{\bold G_{#1,#2}}
\define\hh#1{\Bbb H^{#1}}
\define\cp#1{\Bbb {CP}^{#1}}
\define\qq#1{\Bbb Q^{#1}}
\define\qp#1{\Bbb {QP}^{#1}}

\define\df#1{\frac{\partial\phi}{\partial #1}}
\define\ddf#1{\frac{\partial^2\phi}{\partial #1^2}}
\define\a{\alpha}
\define\be{\beta}
\define\de{\delta}
\define\De{\Delta}

\define\g{\gamma}

\define\n{\nabla}
\define\pa{\partial}
\define\vep{\varepsilon}

\font\aufont=cmr8

\rightheadtext{entire space-like submanifolds}
\leftheadtext{J. Jost and Y.L. Xin}
\topmatter
\title
Some aspects of the global geometry of entire space-like submanifolds \endtitle
\author J. Jost 
\\{\aufont Max-Planck Institute for Mathematics in the Sciences, \\
D-04103 Leipzig, Germany}\\
\\and \\
\\Y.L.Xin
\\{\aufont Institute of Mathematics, Fudan University,\\ Shanghai 200433, China}
\endauthor
\abstract
We prove some Bernstein type 
theorems for entire space-like subma-nifolds in 
pseudo-Euclidean space and as a
corollary, we give a new proof of the Calabi-Pogorelov theorem
for Monge-Amp\`ere equations.
\endabstract
\thanks
1991 Mathematics Subject Classification Primary 53C40, 58E20 \endthanks
\thanks The second-named author thanks the Max Planck Institute 
for Mathematics in the Sciences for its hospitality and support.\endthanks
\endtopmatter

\document
\head 1. Introduction\endhead
\medskip

The search for Bernstein theorems, i.e. theorems stating that, perhaps 
under suitable conditions, entire minimal
graphs or their higher codimensional analogues in Euclidean space,
necessarily are flat and planar, has been a
central topic in geometric analysis and led to many
important insights in the regularity theory of PDEs.
>From that perspective, it might look like a curiosity
to study the analogous question for space-like
entire minimal submanifolds that can be represented
by a graph over a linear subspace in pseudo-Euclidean space. It turns out, 
however, that this situation leads to a rich mathematical structure of its 
own, with many aspects not shared by their Euclidean counterparts. 

The investigations started with a paper of E. Calabi [C1].
Let $ f: {\Bbb R}^m \to {\Bbb R} $ be a smooth function.
If its graph $ (x^1, \dots, x^m ; f (x^1, \dots, x^m)) $ in Minkowski space 
defines a space-like extremal hypersurface, then $ f $ satisfies
$$
 ( 1 - \vert \nabla f \vert^2) 
 \sum_{i=1}^m {\partial^2 f \over \partial {x^i}^2} +
 \sum_{i, j=1}^m {\partial f\over \partial x^i} {\partial f \over \partial x^j} 
 {\partial^2 f\over \partial x^i \partial x^j} = 0, \quad
 \vert \nabla f \vert < 1,
 \tag 1.1
$$

E. Calabi raised the Bernstein type problem and
proved that (1.1) has only linear entire solutions for
$ m \leq 4 $. Several years later the problem was
settled for all $ m $ by Cheng-Yau [C-Y].

The issue of space-like hypersurfaces of
nonzero constant mean curvature becomes more complicated.
The results here depend on conditions about the image under the Gauss map
of the hypersurface [Cho-Tr] [X1] [X-Y].

If we wish to study such entire minimal submanifolds
of higher codimension, we need to study maps
$$
 f : {\Bbb R}^m \to {\Bbb R}^n 
$$
solving the generalization of (1.1) which reads as
$$\sum_{i,j=1}^m g^{ij} {\pa^2 f^s \over \pa x^i\pa x^j} = 0,$$
where $g_{ij} = \de_{ij} - \sum_{s=1}^n f^s_{x^i}f^s_{x^j}$ 
and $\left(g^{ij}\right)$  is the inverse matrix of $\left(g_{ij}\right).$ 
We let $ M $ be the graph of $ f $. A natural approach is to investigate 
the Gauss map
$$
 \gamma : M \to \grs{m}{n}^n,
$$
associating to each point in $ M $ its tangent space, 
considered as an element of the pseudo-Grassmannian $ \grs{m}{n}^n $.
In contrast to its dual and Euclidean counterpart $ \grs{m}{n} $,
$ \grs{m}{n}^n $ is a symmetric space of noncompact type.

Therefore, in particular, it carries a complete Riemannian
metric of nonpositive curvature. This is much
better adapted to the geometry of the Gauss map as a harmonic map
than the (partly) positive curvature of $ \grs{m}{n} $,
and this can be considered as the reason why
one has stronger Bernstein theorems in the
pseudo-Euclidean than in the Euclidean case. Related
to this fact is that for a space-like submanifold
in pseudo-Euclidean space with parallel mean
curvature, its Ricci curvature is bounded from below.
This is a main advantage compered with the Euclidean
space as the ambient space. By Cheng-Yau's
method we have two estimates for the squared norm of
the second fundamental form. One is in terms of
the mean curvature, and the other estimate
is in terms of its mean curvature and the image diameter of its Gauss map. 

In the present situation, however, an entire solution
not necessarily defines a complete manifold as in the ambient Euclidean case. 
The completeness argument thus becomes a key issue. Fortunately, Cheng-Yau's 
method allows an extension to the higher codimensional case. We also obtain a 
gradient estimate for the pseudo-distance on a space-like $m$-submanifold in
pseudo-Euclidean space $ \Bbb R^{m+n}_n $ with index $ n $.

All the geometric conclusions in the paper stem from those estimates. 

We obtain a general Bernstein type theorem stating
that any space-like submanifold in pseudo-Euclidean space
with parallel mean curvature that is closed w.r.t. the
Euclidean topology and whose Gauss map is bounded
is necessarily planar. 

Apparently, Hitchin [H] was the first to observe
the connection between Lagrangian
minimal graphs in pseudo-Euclidean space and Monge-Amp\`ere 
equations. Namely, if $ F : \Bbb R^m \to \Bbb R $ satisfies
the Monge-Amp\`ere equation
$$
\det\left({\partial^2 F \over\partial x^i \partial x^j}\right) = \hbox{ const.}
 \eqno(1.2)
$$
then the graph of its gradient defines a special (i.e. minimal)
Lagrangian submanifold of $ \Bbb R^{2m}_m $. This submanifold
is space-like precisely if $ F $ is convex. Therefore,
we may apply our Bernstein theorem to obtain
a new proof of the famous theorem of Calabi [C2]
(dimension $ \leq 4 $) and Pogorelov [P] (any dimension) that the only entire 
convex solutions are quadratic polynomials which is a
fundamental result in affine differential geometry.
There may exist further connections with Lagrangian
geometry related to the mirror symmetry conjecture.
The starting point is McLeans's construction [M]
of the moduli space $ M $ of special Lagrangian submanifolds
of a Calabi-Yau manifold. McLean constructed a
natural Riemannian metric on this moduli space $ M $.
The key result for us now is Hitchin's [H]
construction of a natural embedding of $ M $ as a
Lagrangian submanifold of pseudo-Euclidean space so that its
space-like metric is precisely McLean's metric.

Now, this embedding of $ M $ in general is not
minimal, but if it is, in view of our results,
this has strong geometric consequences for the
space of special Lagrangian submanifolds of the original
Calabi-Yau manifold. We therefore believe that
exploring these connections between Lagrangian geometry
and Bernstein theorems in pseudo-Euclidean spaces is
a rewarding topic for future research.

\bigskip\bigskip

\head 2. Estimates of the second fundamental form\endhead
\medskip

Let $\ir{m+n}_n$ be an $(m+n)$-dimensional pseudo-Euclidean space
of index $n.$
Let $M$ be a space-like oriented $m$-submanifold in $\ir{m+n}_n.$ 
Choose a local Lorentzian frame field 
$\{e_i,e_s\}$ along $M$ with dual frame field 
$\{\omega_i,\omega_s\}$, such that the $e_i$ are tangent vectors to $M$. 
We agree with the following range of indices
$$A,\,B,\,C,\cdots=1,\cdots,m+n;$$
$$i,\,j,\,k\cdots=1,\cdots, m;\;s,\,t,\cdots=m+1,\cdots,m+n.$$

The induced Riemannian metric of $M$ is given by $ds_M^2 =\sum_i\omega_i^2$
and the induced structure equations of $M$ are
$$ \align
   & d\omega_i = \omega_{ij}\wedge\omega_j,\qquad
                 \omega_{ij}+\omega_{ji} = 0,\cr
   &d\omega_{ij}= \omega_{ik}\wedge\omega_{kj}-\omega_{is}\wedge\omega_{sj},\cr
   &\Omega_{ij} = d\omega_{ij}-\omega_{ik}\wedge\omega_{kj}
                = -\frac 12 R_{ijkl}\omega_k\wedge\omega_l.
\endalign $$
By Cartan's lemma we have
$$\omega_{si} = h_{sij}\omega_j,$$
where the $h_{sij}$ are the components of the second fundamental form of
$M$ in $\ir{m+n}_n$.  The mean curvature vector of $M$ in $\ir{m+n}_n$ is 
defined by
$$H = \frac 1mh_{sii}e_s.$$
We also have the Gauss equation
$$R_{ijkl} = -(h_{sik}h_{sjl} - h_{sil}h_{sjk}),\tag 2.1$$
and the Ricci curvature
$$R_{ij} = R_{kikj} = -( h_{skk}h_{sij} - h_{ski}h_{skj}),$$
from which it follows that 
$$\text{Ric}_M\ge -\frac 14 m^2|H|^2.\tag 2.2$$  
There is an induced connection  on the normal bundle $NM$
in $\ir{m+n}_n$. We have
$$ \align
   &d\omega_{st} = -\omega_{sr}\wedge\omega_{rt} + \Omega_{st},\cr
   &\Omega_{st} = -\frac 12 R_{stij}\, \omega_i\wedge\omega_j,\cr
   &R_{stij} = (h_{ski}h_{tkj} - h_{skj}h_{tki}).\tag 2.3
\endalign $$
The covariant derivative of $h_{sij}$ is given by
$$ h_{sijk}\,\omega_k = dh_{sij} + h_{slj}\,\omega_{li} 
                                 + h_{sil}\,\omega_{lj} 
                                 - h_{tij}\,\omega_{ts}.\tag 2.4$$
It is easily seen that $h_{sijk}=h_{sikj}$, so $h_{sijk}$ is
symmetric in $i, j, k.$ If 
$$DH = \frac 1mh_{siik}\,\omega_k e_s \equiv 0,\tag 2.5$$
then $M$ is called a space-like submanifold with parallel mean curvature.   
In this case we have [X2]

$$ \frac 12\De S\ge \sum h_{sijk}^2
       - m|H|S^{\frac 32} +\frac 1n S^2,\tag 2.6$$
where $S$ is the squared norm of the second fundamental form.

The space-like $m-$planes in $\ir{m+n}_n$ form the pseudo-Grassmannian
$\grs{m}{n}^n.$ It is a symmetric space of noncompact type which is the
noncompact dual space of the Grassmannian manifold $\grs{m}{n}.$ The
canonial Riemannian metric on $\grs{m}{n}^n$ is given by
$$ds_G^2=\sum_{s,i}(\omega_{si})^2.$$

Let 0 be the origin of $\ir{m+n}_n$. Let $SO^0(m+n,n)$ denote the identity 
component of the Lorentzian group $O(m+n,n)$. $SO^0(m+n,n)$ can be viewed as 
the manifold consisting of all Lorentzian frames $(0;e_i,e_s),$ and 
$SO^0(m+n,n)/SO(m)\times SO(n)$ can be viewed as $\grs{m}{n}^n$. Let 
$P=\{(x;e_1,\cdots,e_m);x\in M, e_i\in T_xM\}$ be the principal
bundle of orthonormal tangent frames over $M$,
$Q=\{(x;e_{m+1},\cdots,e_{m+n}); x\in M, e_s\in N_xM\}$ be the principal
bundle of orthonormal normal frames over $M$, then $\bar \pi : 
P\oplus Q \to M$ is the projection with fiber $SO(m)\times SO(n),$
$\italic i:P\oplus Q\hookrightarrow SO^0(m+n,n)$ is the natural inclusion.

We define the generalized Gauss map $\gamma : M \to \grs{m}{n}^n$ by 
$$\gamma (x) = T_xM \in \grs{m}{n}^n$$
via parallel translation in $\ir{m+n}_n$ for $\forall x\in M$. Thus,
the following commutative diagram holds
$$ \CD
   P\oplus Q @>\italic i>> SO^0(m+n,n) \cr
    @V\bar\pi VV  @VV\pi V \cr
   M@>\gamma >> \grs{m}{n}^n
   \endCD $$
With respect to the canonical metric $ds_G^2$ of $\grs{m}{n}^n$ the Levi-Civita
connection is given by
$$ \omega_{(si)(tj)} = \delta_{st}\omega_{ij} -
                       \delta_{ij}\omega_{st}.$$
Using the above diagram, we have
$$ \gamma^*\omega_{si} = h_{sij}\omega_j. \tag 2.7$$

Let $r,\;\tilde r$ be the respective distance functions on $M$ and
$\grs{m}{n}^n$ relative to fixed points 
$x_0 \in M,\;\tilde x_0 \in \grs{m}{n}^n.$
Let $B_a$ be a closed geodesic ball of radius $a$ and centered at $x_0.$ 
Define the maximum modulus of the Gauss map 
$\gamma: M \to \grs{m}{n}^n$ on $B_a$ 
by
$$\mu(\gamma,a)\mathop{=}\limits^{def.}\max\,\{\tilde r(\gamma(x));\,
x \in B_a \subset M \}.$$
For a fixed positive number $a,$ choose $b\ge\mu^2(\gamma,a).$ Define $f:
B_a\to \ir{}$ by
$$f=\frac{(a^2-r^2)^2S}{(b-h\circ\gamma)^2},$$
where $h=\tilde r^2.$ By applying the maximum principle to $f$ we can derive an 
estimate for $S$, the squared norm of the second fundamental form of $M$ 
in $\ir{m+n}_n,$ in terms of the mean curvature and the image diameter of the 
Gauss map. In fact, we have [X2] that for any $x \in B_a$
$$S(x)\le k\Bigl(\,\frac {(8\mu\,a+m\,a^2|H|)^2\mu^4}
                        {(2+\frac 1n\mu^2)^2(a^2-r^2)^2}
      +\frac {(2(m+4)a^2+m(m-1)|H|\,a^3)\mu^2}
              {(2+\frac 1n\,\mu^2)(a^2-r^2)^2}\,\Bigr),\tag 2.8$$
where $k$ is an absolute constant.  In what follows $k$ may
be different in different inequalities. 

Consider the auxiliary function
$$f = (a^2-r^2)^2S$$
on a geodesic ball $B_a$ of radius $a$ and centered at $x_0\in M.$ By a 
similar method we can obtain an estimate in terms only of the mean curvature
of $M$ in $\ir{m+n}_n:$
$$S(x)
\le k\,\frac {m^2n^2|H|^2a^4+mn(m-1)|H|\,a^3+2n(m+4)\,a^2}{(a^2-r^2)^2},
\tag 2.9$$
for all $x\in B_a\subset M$.

\bigskip\bigskip
\head 3. Completeness \endhead
\medskip
In this section we generalize the argument of Cheng-Yau [C-Y] to higher
codimension  with some technical modifications. 
 
Let $M$ be a space-like submanifold in pseudo-Euclidean space $\ir{m+n}_n$
with index $n.$
Let $X=(x_1, \cdots, x_m;\, y_1, \cdots, y_n)$ be the position vector of $M.$
Define the pseudo-distance function on $M$ by $\left <X, X\right >
=\sum_i x_i^2-\sum_s y_s^2.$ It is non-negative because $M$ is space-like.

\proclaim{Proposition 3.1}
If $M$ is closed with respect to the Euclidean topology, then when $0\in M,$ 
$z=\left <X,X\right >$ is a proper function on $M.$
\endproclaim

\demo{proof}
Let $\bar c= \inf\{ c;\;\text{the set where}\;\left <X, X\right >\le c 
\;\text{is compact} \}.$ Then we will show $\bar c=\infty.$

Let $\ir{2}_1\subset \ir{m+n}_n$ be a Minkowski plane. Since 
$0\in \ir{2}_1\cap M$ and $M$ is space-like, $M$ meets $\ir{2}_1$ transversally.
It follows that there are positive constants $\vep_1, \vep_2$ and $\vep_3$
such that for $(x, y)\in \ir{2}_1\cap M$ and $\sum_i x_i^2 = \vep_1,$ we have
$\vep_2 \ge\left <X, X\right >\ge \vep_3.$

Suppose $\bar c < \infty.$ By the assumption that $M$ is closed with respect to
the Euclidean topology we have a sequence of points 
$(x_1^\a, \cdots, x_n^\a; y_1^\a, \cdots, y_m^\a)$ in $M$ such that 
$$\align
\sum_i(x_i^\a)^2&\to\infty,\cr
\sum_s(y_s^\a)^2&\to\infty,\cr
\sum_i(x_i^\a)^2&-\sum_s(y_s^\a)^2<\bar c.
\endalign$$

Choose $\a$ sufficiently large, such that

$$\sqrt{\sum_i(x_i^\a)^2} 
   > \vep_1^{\frac 12}\vep_3^{-1}(\bar c+2\vep_2).\tag3.1$$
By an action of $SO(m)\times SO(n)$ we have new coordinates of $\ir{m+n}_n$
such that the point $(x_1^\a, \cdots, x_m^\a; y_1^\a, \cdots, y_n^\a)$ becomes
$(\sqrt{\sum_i(x_i^\a)^2}, 0\cdots,0; \sqrt{\sum_s(y_s^\a)^2}, 0, \cdots, 0)$ in 
the new coordinates. For simplicity it is denoted by $(x_1^\a, 0, \cdots, 0; 
y_1^\a, 0, \cdots, 0)$ with $y_1^\a>0$ and (3.1) becomes
$$x_1^\a>\vep_1^{\frac 12}\vep_3^{-1}(\bar c+2\vep_2).\tag3.2$$

Let $P^\a$ be the Minkowski $2-$plane spanned by the $x_1-$axis and the
$y_1-$axis. By the previous argument $P^\a$ intersects $M$ in a point 
$(x_1^0, 0, \cdots, 0; y_1^0, 0,\cdots, 0)$ with $(x_1^0)^2= \vep_1$ and
$\vep_2\ge (x_1^0)^2-(y_1^0)^2\ge \vep_3.$

Since $M$ is space-like, the point $(x_1^\a, 0, \cdots, 0; y_1^\a, 0,\cdots, 0)$
can not lie in the light cone of $(x_1^0, 0, \cdots, 0; y_1^0, 0, \cdots, 0).$
Therefore,
$$\bar c+(x_1^0)^2-(y_1^0)^2
    \ge 2x_1^\a(x_1^0-y_1^0)+2(x_1^\a-y_1^\a)y_1^0,\tag 3.3$$
and
$$\bar c (x_1^\a+y_1^\a)^{-1} \ge x_1^\a -y_1^\a \ge 0. \tag 3.4$$
>From $x_1^0=\sqrt {\vep_1} > |y_1^0|$ and $(x_1^0)^2-(y_1^0)^2 \ge \vep_3$
we have
$$x_1^0-y_1^0 > 2^{-1}\vep_1^{-\frac 12}\vep_3. \tag 3.5$$
Substituting (3.4) and (3.5) into (3.3) gives
$$\bar c +\vep_2 \ge 
\vep_1^{-\frac 12}\vep_3 x_1^\a - \bar c\;\vep_1^{\frac 12}(x_1^\a)^{-1},$$
namely,
$$\vep_1^{-\frac 12}\vep_3(x_1^\a)^2 - (\bar c+\vep_2) x_1^\a
-\bar c\;\vep_1^{\frac 12} \le 0.$$
It follows that
$$x_1^\a < \vep_1^{\frac 12}\vep_3^{-1}(\bar c + 2\vep_2),$$
which contradicts (3.2) and the proof is complete.\hfill$\square$

\enddemo

Now let us study $X: M \to \ir{m+n}_n$ being a space-like submanifold with
parallel mean cuvature. Choose a Lorentzian frame field $\{e_i, e_s\}$ along
$M,$ such that $e_i$ are tangent vectors to $M$ with $\n_{e_i}e_j =0$ at the 
considered point. We need also carefully choose the normal vectors. This is
the main technical point to generalize Cheng-Yau's proof to higher codimension.

Let $\bar X= X-\left<X, e_i\right>e_i.$ At a point, say $q$, choose
$$e_{m+1}=\frac {\bar X}{|\bar X|},$$
then choose other normal vectors $e_{m+2}, \cdots, e_{m+n}$, so that they are
all mutually orthogonal and then expand
them around the point $q$ to form a local normal frame field.

Let $z=\left <X, X\right>$ be the pseudo-distance function on $M.$ 
Then
$$z_i\mathop{=}\limits^{def.} e_i(X)= 2\left<X, e_i\right>,\tag 3.6$$
$$z_{ij}\mathop{=}\limits^{def.} Hess (z)(e_i,e_j)
=2 \left(\de_{ij} - \left<X, e_s\right> h_{sij}\right),\tag 3.7$$
$$\De z = 2 m - 2 m \left<X, e_s\right>H_s.\tag 3.8$$

On the compact set $\{z\le k\}$ in $M$ for some $k$ define a function
$$f =\frac {|\n z|^2}{(z+1)^2}\exp\left(\frac{-c}{k-z}\right),\tag 3.9$$
where $c$ will be chosen later. It attains its maximum at a point $q.$ Then,
$$\n f (q) = 0,$$
$$\De f (q) \le 0.$$
By computations we have that at $q$
$$\align
&2 z_jz_{ij} - g |\n z|^2 z_i = 0,\tag 3.10\cr
&2\sum_{ij}z_{ij}^2 + 2 z_jz_{iji} - g' |\n z|^4 
   - g ( 2 z_{ij}z_i z_j +|\n z|^2 \De z ) \le 0,\tag 3.11
\endalign$$   
where
$$g = \frac 2{z+1} + \frac c{(k-z)^2}.$$
By the choice of the normal vectors $e_s$ at the point $q$ (3.7) reduces to
$$z_{ij}=2\de_{ij} - 2\left<X, e_{m+1}\right> h_{m+1\,ij}\tag 3.12$$
It follows that
$$h_{m+1\, ij}=\frac {2\de_{ij}-z_{ij}}{2\left<X, e_{m+1}\right>}, \tag 3.13$$
$$\De z = 2m - 2 m\left<X, e_{m+1}\right>H_{m+1}, \tag 3.14$$
and (3.6) means
$$z=\frac 14|\n z|^2 - \left<X, e_{m+1}\right>^2. \tag 3.15$$

By Schwarz inequlity (see Lemma 2 in [Y]), we have 
$$\sum_{ij}z_{ij}^2\ge
 \frac {2m-1}{2m-2}\sum_i\left(\sum_jz_{ij}z_j\right)^2|\n z|^{-2}
     -\frac 1{m-1}(\De z)^2. \tag 3.16$$  
Substituting (3.10), (3.14) and (3.15) into (3.16) yields
$$\sum_{ij}z_{ij}^2\ge
 \frac {2m-1}{8m-8} g^2|\n z|^4
   -\frac 1{m-1}\left(2m-mH_{m+1}\left(|\n z|^2-4 z\right)^{\frac 12}\right)^2. 
\tag 3.17$$
Noting the Ricci formula, the Gauss equation (2.1) and
$$(\De z)_i =m h_{sij}z_jH_s,$$
we have
$$z_j z_{iji} = h_{ski} h_{skj} z_i z_j \ge h_{m+1\, ki} h_{m+1\,kj} z_i z_j.$$
At the point $q$ we can use (3.13), and the above expression becomes
$$z_j z_{iji}
  \ge 4 - 2 g |\n z|^2 + \frac 14 g^2 |\n z|^4. \tag 3.18$$ 
Substituting (3.10), (3.13), (3.14), (3.17) and (3.18) into (3.11), we have
$$\aligned
0 \ge
\left(\frac 1{4(m-1)}g^2-g'\right)|\n z|^4
&-g\left (4+2m+m|H_{m+1}|(|\n z|^2-4z)^{\frac 12}\right)|\n z|^2\cr
 &+ 8 -\frac 2{m-1}\left (2 m+m|H_{m+1}|(|\n z|^2-4z)^{\frac 12}\right)^2. 
\endaligned\tag 3.19$$
The coefficient of $|\n z|^4$ is
$$\frac 1{4(m-1)}
\left(\frac 4{(z+1)^2}+\frac {c^2}{(k-z)^4}+\frac {4c}{(z+1)(k-z)^2}\right)
                                +\frac 2{(z+1)^2}-\frac {2c}{(k-z)^3}.
\tag 3.20$$
Choose  $ c = 8(m-1)k$
and so that
$$ (3.20) \ge \frac 2{(z+1)^2}.$$
Hence, at the point $q$
$$\frac 2{(z+1)^2}|\n z|^4 \le
g (4+2m +m |H||\n z|) |\n z|^2 + \frac 2{m-1}( 2 m + m|H||\n z|)^2. $$
This means that
$$\align
f^2\le \frac {16m^2}{(m-1)(z+1)^2}\exp\left(\frac {-2c}{k-z}\right)&
+\frac {16m^2}{(m-1)(z+1)}\exp\left(\frac {-3c}{2(k-z)}\right)|H|f^{\frac 12}\cr
  +4(m+2)g\exp\left(\frac {-c}{k-z}\right) f 
    &+\frac {4m^2}{m-1}\exp\left(\frac {-c}{k-z}\right) |H|^2 f \cr
        +2\,m\,g\, (z+1)&\exp \left(\frac {-c}{2(k-z)}\right) |H| f^{\frac 32}.
\endalign$$
We then can find a constant $P$ depending only on $m$ so that
$$f^2\le \frac 14 P\left((z+1)^{-2}+|H|(z+1)^{-1} f^{\frac 12}
                   +(1+|H|^2) f+|H| f^{\frac 32}\right).$$
It follows that
$$\sup_{z\le k} f\le \max\left\{P^{\frac 12}\sup_{z\le k}(z+1)^{-1},\,
                   P^{\frac 23} |H|^{\frac 23}\sup_{z\le k}(z+1)^{- \frac 23},\,
                              P(1+|H|^2),\, P^2|H|^2\right\}.$$
Now, we arrive at the following conclusion.
\proclaim{Proposition 3.2}
Let $M$ be a space-like submanifold in pseudo-Euclidean space $\ir{m+n}_n$                              
of index $n$ with parallel mean cuvature. Let $z$ be the pseudo-distance 
function on $M.$ If for some $k>0,$ the set $\{z\le k\}$ is compact, then 
there is constant $b$ depending only on the dimension $m$ and the norm of 
the mean curvature $|H|,$  such that for all $x \in M$ with $z(x)\le \frac k2,$
$$|\n z| \le b (z+1).\tag 3.21$$\hfill$\square$
\endproclaim

Without loss of  generality we assume that $0 \in M.$ If $M$ is closed 
with respect to the Euclidean topology, then $z$ is a proper function on $M$ by
Proposition 3.1 and (3.21) is valid for any $k$. Let $\g: [0, r] \to M$ be a 
geodesic on $M$ issueing from the origin $0.$ Integrating (3.21) gives
$$z(\g(r))+1 \le \exp(br),$$
which forces $M$ to be complete. In summary we have
\proclaim{Theorem 3.3}
Let $M$ be a space-like submanifold in the pseudo-Euclidean space $\ir{m+n}_n.$
Assume that $M$ is closed with respect to the Euclidean topology and its mean 
cuvature is parallel. Then $M$ is complete with respect to the induced metric 
from the ambient space.
\endproclaim\hfill$\square$

\bigskip\bigskip

\head 4. Bernstein type theorems \endhead
\medskip

We are now in a position to prove some theorems.

\proclaim{Theorem 4.1}
Let $M$ be a space-like $m$-submanifold in pseudo-Euclidean space $\ir{m+n}_n$
with index $n.$ Assume that
\roster
\item $M$ is closed with respect to the Euclidean topology;
\item $M$ has parallel mean curvature;
\item the image under the Gauss map from $M$ into $\grs{m}{n}^n$ is bounded.
\endroster
Then $M$ has to be a linear subspace.
\endproclaim

\demo{Proof}
Let $\tilde x_0\in \grs{m}{n}^n$ and $R$ be a positive number
and large enough such that the image under the Gauss map $\g: M\to \grs{m}{n}^n$
is contained in the geodesic ball $B_R(\tilde x_0)$. Since the mean curvature
is parallel, the Gauss map is harmonic, and so we have a harmonic
map $\g: M\to B_R(\tilde x_0)\subset \grs{m}{n}^n.$ On the other hand,
from (2.2) we know that the Ricci curvature of $M$ is bounded from below.
So, we can use the maximum principle to conclude that (see Thm 3.10 in [X3])
for the energy density of $\g,$
$$\inf e(\g)=0.$$ 
>From (2.7)
$$\inf S=0.\tag 4.1$$
But, on the other hand, by the Schwarz inequality and the assumption of 
parallel mean curvature
$$const. =|H|^2\le \frac 1m S.\tag 4.2$$
(4.1) and (4.2) force that $H= 0.$ Now, we use the estimate
(2.8) and obtain

$$S(x)\le k\left(\frac{32 a^2 R^6 }{(2+\frac 1n R^2)^2(a^2-r^2)^2}
+\frac{(m+4)a^2 R^2}{(2+\frac 1n R^2)(a^2-r^2)^2}\right).\tag 4.3$$
It is valid on a geodesic ball $B_a(x_0)\subset M.$
By Theorem 3.3 $M$ is complete and we can fix $x$ and let $a$ tend to 
infinity in (4.3). Thus, $S(x)=0$ for all $x\in M.$ The proof is complete.
\hfill$\square$\enddemo

Now, we study the following special case. Let $\ir{2m}_m$ be the 
pseudo-Euclidean $2m$ space with index $m.$ Let 
$(x, y) = (x^1, \cdots, x^m;\, y^1, \cdots, y^m)$
be  null coordinates; this means that the indefinite metric is defined by
$$ds^2=\frac 12\sum_i dx^i dy^i.\tag 4.4$$

Let $F$ be a smooth convex function. We consider the graph $M$ of $\n F,$ 
defined by 
$$(x^1, \cdots, x^m;\,\frac {\pa F}{\pa x^1}, \cdots, \frac {\pa F}{\pa x^m}).$$
The induced Rimannian metric on $M$ is defined by
$$ds^2 = \frac {\pa^2 F}{\pa x^i\pa x^j} dx^i dx^j.\tag 4.5$$ 
The underlying Euclidean space $\ir{2m}$ of $\ir{2m}_m$ has the usual complex 
structure. It is easily seen that $M$ is a Lagrangian submanifold in $\ir{2m}.$
Let us derive the condition on $F$ for $M$ being an extremal submanifold
in $\ir{2m}_m.$

Choose a tangent frame field  $(e_1, \cdots, e_m)$ along $M,$ where
$$e_i = 
\frac {\pa}{\pa x^i}+\frac {\pa^2F}{\pa x^i\pa x^j}\frac {\pa}{\pa y^j}.$$
Obviously,
$$\left<e_i, e_j\right>=\frac {\pa^2 F}{\pa x^i\pa x^j}.$$ Let 
$(n_i, \cdots, n_m)$ be the normal frame field of $M$ in $\ir{2m}$ defined by
$$n_i=\frac {\pa}{\pa x^i}-\frac {\pa^2F}{\pa x^i\pa x^j}\frac {\pa}{\pa y^j}$$
with
$$\left<n_i, n_j\right>= - \frac {\pa^2 F}{\pa x^i\pa x^j}.$$
Thus, $M$ is space-like precisely if $F$ is convex.
By direct computations
$$\align\n_{e_i}e_j&=
\n_{\frac {\pa}{\pa x^i}+\frac {\pa^2F}{\pa x^i\pa x^k}\frac {\pa}{\pa y^k}}
\left(\frac {\pa}{\pa x^i}+\frac {\pa^2F}{\pa x^i\pa x^l}\frac {\pa}{\pa y^l}
\right)\cr
&= \frac {\pa^3F}{\pa x^i\pa x^j\pa x^l}\frac {\pa}{\pa y^l}\cr
&=\frac 12\frac {\pa^3F}{\pa x^i\pa x^j\pa x^l}\;g^{lk}\;e_k
  -\frac 12\frac {\pa^3F}{\pa x^i\pa x^j\pa x^l}\;g^{lk}\;n_k,
  \endalign$$
where $g^{ij}$ denotes the elements of the inverse matrix of 
$(g_{ij}) \mathop{=}\limits^{def.}\left(\frac {\pa^2 F}{\pa x^i\pa x^j}\right)$. 
It follows that the second fundamental form of $M$ in $\ir{2m}_m$
$$B_{ij} \mathop{=}\limits^{def.}\left(\n_{e_i}e_j\right)^{\perp}
  = -\frac 12\frac {\pa^3F}{\pa x^i\pa x^j\pa x^l}\;g^{lk}\;n_k, \tag 4.6$$ 
and the mean curvature vector
$$H \mathop{=}\limits^{def.}\frac 1m\sum_iB_{ii}
 = - \frac 1{2m g}\;\frac {\pa g}{\pa x^l}\;g^{lk}\;n_k, \tag 4.7$$
 where $g = \det (g_{ij}).$

\proclaim{Theorem 4.2} Let $M$ be a space-like extremal $m-$submanifold 
in $\ir{m+n}_n.$ If $M$ is closed with respect to the Euclidean topology, then 
$M$ has to be a linear subspace. In particular, when such an $M$ is defined by 
the graph  $(x; \n F)$ of the gradient $\n F$ of a smooth function 
$F: \ir{m}\to \ir{}$
in null coordinates $(x; y)$ in $\ir{2m}_m,$ then $F$ has to be a quadratic 
polynomial.
\endproclaim

\demo{Proof}
By Theorem 3.3 $M$ is complete. On the other hand, substituting $H=0$ into
(2.9) yields
$$S(x)=k\frac {2n(m+4)\,a^2}{(a^2-r^2)^2}.\tag 4.8$$
We fix $x$ and let $a$ go to infinity in (4.8). Hence, $S(x)=0$ for any
$x\in M.$ We complete the proof of the theorem.\hfill$\square$
\enddemo

(4.7) reveals that $H=0$ is equivalent to the Monge - Amp\`ere equation
$$\det\left(\frac {\pa^2 F}{\pa x^i\pa x^j}\right) = \text{const.}\;.\tag 4.9$$ 
We thus obtain an alternative proof of the following famous result  
shown by Calabi (for $m\le 5$) [C2] and Pogorelov [P] (for all dimensions).

\proclaim{Corollery 4.3}
The only entire convex solutions to (4.9) are  quadratic polynomials.
\endproclaim

\remark{Remarks}
1. In [A] and [I], some results were proved for complete space-like submanifolds.
It was shown in [A] that completeness implies that the manifolds are closed with
respect to the Euclidean topology. Therefore, the present results are 
generalizations of their results.

2. A convex solution to the Monge - Amp\`ere eqation (4.9) represents
an improper affine hypersphere in affine differential geometry. The present 
discussion shows the close relationship between Lagrangian extremal submanifolds 
and affine hypersurfaces.
\endremark

\bigskip\bigskip

\head 5. Final remarks\endhead
\medskip
Let $L$ be a compact special Lagrangian submanifold of a Calabi - Yau manifold 
$Y.$ Let $M$ be the moduli space of the special Lagrangian submanifolds near 
$L.$ McLean [M] showed that $M$ is a smooth manifold of dimension $\be (L)
= \dim (H^1(L, \ir{})).$ He also defined a natural Riemannian metric on $M.$
Hitchin [H] then studied this moduli space $M.$  He showed that there is 
natural embedding of the moduli space $M$ as a Lagrangian submanifold in
the product $H^1(L,\ir{})\times H^{m-1}(L,\ir{})$ (where $m=\dim L$) of
two dual vector spaces and that McLean's metric is the metric induced by the
ambient pseudo-Euclidian metric. He also showed that as a Lagrangian
submanifold $M$ is defined locally by graph of the gradient of a function
$F.$  So, we are in the situation studied in the previous section.

Therefore, the curvature properties of the moduli space $M$ can be 
determined by our previous calculations. From the Gauss equation and (4.6)
we obtain the Riemannian curvature, the Ricci curvature and the scalar curvature
of the moduli space $M$ with respect to McLean's metric as follows.
$$R_{ijkl}
= - \frac 14\;g^{st}\frac {\pa^3 F}{\pa x^s\pa x^i\pa x^k}
                    \frac {\pa^3 F}{\pa x^t\pa x^j\pa x^l}
  + \frac 14\;g^{st}\frac {\pa^3 F}{\pa x^s\pa x^i\pa x^l}
                    \frac {\pa^3 F}{\pa x^t\pa x^j\pa x^k} ,$$
                    
$$R_{ik}
= - \frac 1{4\;g}\;g^{st}\frac {\pa^3 F}{\pa x^s\pa x^i\pa x^k}
                         \frac {\pa g}{\pa x^t}
  +  \frac 14\;g^{st}\;g^{jl} \frac {\pa^3 F}{\pa x^s\pa x^i\pa x^l}                              
                              \frac {\pa^3 F}{\pa x^t\pa x^j\pa x^k} ,$$
and
$$R = - \frac 14\;g^{st}\;\frac {\pa\ln g}{\pa x^s}\,\frac {\pa\ln g}{\pa x^t}
 +  \frac 14\;g^{st}\,g^{jl}\,g^{ik}\frac {\pa^3 F}{\pa x^s\pa x^i\pa x^l} 
                                    \frac {\pa^3 F}{\pa x^t\pa x^j\pa x^k}.$$
                                    
It is interesting to observe that when the moduli space $M$ not only is  
Lagrangian, but also is special, in this case by (4.9)                                    
the Ricci curvature of the moduli space is nonnegative.


\vskip0.5in

\Refs
\medskip
\widestnumber\key{Cho-Tr}

\ref\key A
\by Aiyama, ~R.~
\paper The generalized Gauss map of a space-like submanifold with parallel
mean curvature vector in pseudo-Euclidean space
\jour Japan J. Math.
\vol 20(1)
\yr 1994
\pages 93-114
\endref
\medskip

\ref\key C1
\by Calabi,~E.~
\paper Examples of Bernstein problems for some nonlinear equations
\jour Proc. Symp. Global Analysis U.C.Berkeley
\yr 1968
\endref
\medskip

\ref\key C2
\by Calabi, ~E.~
\paper Improper affine hyperspheres of convex type and generalization
of a theorem by K. J\"ogens'
\jour Michigan Math. J.
\vol 5
\yr 1958
\pages 105-126
\endref
\medskip

\ref\key C-Y
\by Cheng,~S.~Y.~ and Yau,~S.~T.~
\paper Maximal spacelike hypersurfaces in the Lorentz-Minkowski space
\jour Ann. Math.
\vol 104
\yr 1976
\pages 407-419
\endref
\medskip


\ref\key Cho-Tr
\by Choi,~H.~In~and Treibergs,~A.~
\paper Gauss maps of spacelike constant mean curvature
hypersurfaces of Minkowski space
\jour J. Differential Geometry
\vol 32
\yr 1990
\pages 775-817
\endref
\medskip

\ref\key H
\by Hitchin, ~N.~J.~
\paper The moduli space of special Lagrangian submanifolds
\jour Ann. Scuola Norm. Sup. Pisa Cl. Sci 
\vol 25(4)
\yr 1998
\pages 503-515
\endref
\medskip

\ref\key I
\by Ishihara,~T.~
\paper The harmonic Gauss maps in a generalized sense
\jour J. London Math. Soc.
\vol 26, 104-112
\yr 1982
\endref
\medskip

\ref\key M
\by McLean, ~R.~C.~
\paper Deformations of calibrated submanifolds
\jour Comm. Anal. Geom.
\vol 6
\yr 1998
\pages 705-747
\endref
\medskip

\ref\key P
\by Pogorelov, ~A.~V.~
\paper On the improper convex affine hyperspheres
\jour Geometriae Dedicata
\vol 1
\yr 1972
\pages 33-46
\endref 
\medskip

\ref\key X1
\by Xin,~Y.~L.~ 
\paper On Gauss image of a spacelike hypersurface with constant
  mean curvature in Minkowski space
\jour Comment. Math. Helv.
\vol 66, 
\pages 590-598
\yr 1991
\endref
\medskip

\ref\key X2
\by Xin,~Y.~L.~
\paper A rigidity theorem for a space-like graph of higher codimension
\jour manuscripta math.
\vol 103(2)
\yr 2000
\pages 191-202
\endref
\medskip

\ref\key X3
\by Xin,~Y.~L.~
\book Geometry of harmonic maps
\publ Birkh\"auser, PNLDE \; vol. 23
\yr 1996
\endref
\medskip

\ref\key X-Y
\by Xin, ~Y.~L.~and Ye,~Rugang
\paper Bernstein-type theorems for space-like surfaces with parallel
mean curvature
\jour J. reine angew. Math.
\vol 489
\yr 1997
\pages 189-198
\endref
\medskip

\ref\key Y
\by Yau, ~S.~T.~
\paper Harmonic function on complete Riemannian manifold
\jour Comm. Pure Appl. Math.
\vol 28
\yr 1975
\pages 201-228
\endref

\endRefs
\enddocument